\documentclass[a4paper]{article}

\usepackage{graphics,amsmath,amsfonts}

\newtheorem{Theorem}{Theorem}[section]
\newtheorem{Remark}{Remark}[section]

\newcommand{{\cT}}{{\cal T}}
\newcommand{{\cC}}{{\cal C}}

\newcommand{\R}{{\mathbb R}}

\newcommand{\disp}{\displaystyle}

\newcommand{\w}{\widetilde v}

\title{Initialization of the Shooting Method via the Hamilton-Jacobi-Bellman Approach\footnote{Acknowledgments: the authors wish to thank Hasnaa Zidani for proposing the main idea of the paper and for her suggestions.}}

\author{E. Cristiani\footnote{Emiliano Cristiani (corresponding author): CEMSAC, Universit\`a di Salerno, Fisciano (SA), Italy and IAC-CNR, Rome, Italy. E-mail: \texttt{emiliano.cristiani@gmail.com}}, P. Martinon\footnote{Pierre Martinon: INRIA and CMAP \'Ecole Polytechnique, Palaiseau, France. E-mail: \texttt{martinon@cmap.polytechnique.fr}}}

\date{\today}

\begin{document}

\maketitle

\begin{abstract}
The aim of this paper is to investigate from the numerical point of view the possibility of coupling the Hamilton-Jacobi-Bellman
(HJB) equation and Pontryagin's Minimum Principle (PMP) to solve some control problems. 
A rough approximation of the value function computed by the HJB method is used to obtain 
an initial guess for the PMP method. 
The advantage of our approach over other initialization techniques (such as continuation or direct methods) is to provide an initial guess close to the global minimum.
Numerical tests involving multiple minima, discontinuous control, singular arcs and state 
constraints are considered. The CPU time for the proposed method is less than four minutes up to 
dimension four, without code parallelization. 
\end{abstract}
%
%
\section{Introduction}

In optimal control, the so-called direct methods, based on discretization
and nonlinear programming, are currently the most popular.
The development of many powerful codes in the recent years, such as
NUDOCCCS \cite{NUDOCCCS}, MUSCOD \cite{MUSCOD}, or IPOPT \cite{IPOPT}, allowed to solve
difficult and complex problems \cite{CBG02, MP08}.
On the other hand, the indirect methods, based on Pontryagin's Minimum
Principle (PMP), are both fast and accurate but tend to suffer from a great
sensitivity to the initialization.
The aim of this paper is to show that this difficulty can be overcome by
coupling the indirect methods with the Hamilton-Jacobi-Bellman (HJB) approach.

The HJB theory and PMP are usually considered two separate worlds although they deal with the same kind of problems. The theoretical connections between the two approaches are well known \cite{ClVi87,CaFra91,CFS00,CF05}, 
but coupled usage of the two techniques is not common and not completely explored.\\
In this paper we deal with the following controlled dynamics
\begin{equation}\label{dynamicalsystem}
\left\{
\begin{array}{ll}
\dot y(t)=f(y(t),u(t)), & t>0 \\
y(0)=x, & x\in\R^d
\end{array}
\right.
\end{equation}   
where the control variable $u(\cdot)\in\mathcal U:=\{u:\R^+\rightarrow U,~ u \textrm{ measurable}\}$ and $U\subset\R^m$ 
($m\geq 1$). We will denote by $y_x(t;u)$ the solution of the system (\ref{dynamicalsystem}) starting from the point $x$ with control $u$. 
Let $\mathcal C\subset \R^d$ be a given \emph{target}. For any given control $u$ we denote by $t_f(x,u)$ the first time the trajectory $y_x(t;u)$ hits $\mathcal C$ (we set $t_f(x,u)=+\infty$ if the trajectory never hits the target). We also define a \emph{cost functional} $J$ as
\begin{equation}\label{cout}
J(x,u):=\int_0^{t_f(x,u)}\ell(y_x(t;u),u(t))dt\,,
\end{equation}
where $\ell$ is a suitable cost function. 
The final goal is to find
\begin{equation}\label{mainproblem0}
u^*\in\mathcal U \textrm{ such that } J(x,u^*)=\min\limits_{u\in\mathcal U} J(x,u)\,,
\end{equation}
and compute the associated optimal trajectory $y^*_x(t;u^*)$. \\
Finally, we define the \emph{value function} $$\cT(x):=J(x,u^*)\,,\quad x\in \R^d.$$
Choosing $\ell\equiv 1$ in (\ref{cout}) we obtain the classical \emph{minimum time} problem.

\medskip

The PMP approach consists in finding a trajectory which satisfies some necessary conditions. This is done in practice by searching a zero of a certain shooting function, typically with a \mbox{(quasi-)Newton} method. This method is well known and it is used in many applications, see \cite{Pesch, Pesch_ref, Deu04} and references therein. 
The main advantages of this approach lie in its accuracy and in its low numerical complexity. It is worth to recall that the dimension of the nonlinear system for the shooting method is usually $2d$, where $d$ is the state dimension. This is in practice quite low for this kind of problem, therefore fast convergence is expected in case of success, especially if the initial guess is close to the right value. Unfortunately, finding a suitable initial guess can be extremely difficult in practice. 
The algorithm may either not converge at all, or converge to a local minimum of the cost functional.

\medskip

The HJB approach is based on the Dynamic Programming Principle \cite{Bellman}. It consists in characterizing the value function $\cT$ as the solution of a first-order nonlinear partial differential equation. Once an approximation of the value function is computed, they are easily obtained both the optimal control $u^*$ in feedback form and, by a direct integration, the optimal trajectories for any starting point $x\in\R^d$ \cite{F_appendix, F_games}. The method is greatly advantageous because it is able to reach the global minimum of the cost functional, even if the problem is not convex.
The HJB approach allows also to have a global overview of the optimal trajectories and of the reachable sets (or capture basins) i.e. the sets of the points from which it is possible to reach the target in any given time. 
Beside all the advantages listed above, the HJB approach suffers from the well known ''curse of dimensionality'', so in general it is restricted to problems in low dimension ($d\leq 3$).

\medskip

In this paper we couple the two methods in such a way we can preserve the respective advantages. The idea is to solve the problem via the HJB method on a coarse grid, to have in short time a first approximation of the value function and the structure of the optimal trajectory. 
Then, we use this information to initialize the PMP method and compute a precise approximation of the global minimum. To our knowledge this is the first attempt to exploit the connection between HJB and PMP theories from the numerical point of view.

Compared to the use of continuation techniques or direct methods to obtain an estimate of the initial costate, the main advantage of the approach presented here is that the HJB method provides an initial guess close to the global minimum. The main limitation is the restriction with respect to the dimension of the state.  

We consider some known control problems with different specific difficulties: several local minima, several global minima, discontinuous control, presence of singular arcs, and state constraints. In all these problems, we show that combining PMP method with HJB approach leads to a very efficient algorithm.

\section{Preliminaries}
Consider optimal control problems in the general Bolza form, autonomous case, with a fixed or free final time.
$$
\label{mainproblem}
(P) \left \lbrace
\begin{array}{lr}
\min\ J(x,u) =\int_{0}^{t_f(x,u)} \ell(y(t),u(t)) dt   & \mbox{Objective}\\
\dot{y}(t) = f(y(t),u(t))                               & \mbox{Dynamics}\\
u (t)\in U\quad \mbox{for a.e. } t\in (0,t_f(x,u)) & \mbox{Admissible\ Controls}\\
y(0) = x                                     & \mbox{Initial\ Conditions}\\
y(t_f(x,u)) \in \mathcal{C}                    & \mbox{Terminal\ Conditions}\\
\end{array}
\right .
$$
Here $U$ is a compact set of $\R^m$ and the following classical assumptions are satisfied:
\begin{itemize}
 \item[-] $f:\R^d\times U\rightarrow\R^d$ and $\ell:\R^d\times U\rightarrow\R$ are continuous, and are of class $C^1$ with respect to the 
  first variable.
 \item[-] $\cC$ is a closed subset of $\R^d$ for which the property ''a vector is normal to $\mathcal{C}$ at a point of $\mathcal{C}$'' makes sense. For instance, $\cC$ can be described by a finite set of equalities $\{c_i(x)=0\}_i$ or inequalities $\{c_i(x)\le 0\}_i$, with $c_i$ being of class $C^1$ for every $i$, and the classical constraint qualification assumptions hold. For numerical purposes we assume $\cC$ is bounded.
\end{itemize}

\subsection{Pontryagin's Minimum Principle approach}
\label{methodPMP}
We give here a brief overview of the so-called indirect methods for optimal control problems  
\cite{Pontryagin, BryHo75, Pes89}.
We introduce the costate $p$, of same dimension $d$ as the state $x$, and define the Hamiltonian
$$
H(y,p,u,p_0) = p_0 \ell(y,u) + <p,f(y,u)>.
$$

Under the assumptions on $f$ and $\ell$ introduced above, Pontryagin's Minimum Principle states that if 
$(y_x^*,u^*,t_f^*)$ is a solution of ($P$), then there exists $(p_0,p^*) \neq 0$ absolutely continuous such that
\begin{subequations}\label{PMP}
\begin{eqnarray}
& & \dot{y}^*(t) =  H_p(y_x^*(t),p^*(t),u^*(t),p_0),\quad y_x^*(0)=x,\\
& & \dot{p}^*(t) = - H_y(y_x^*(t),p^*(t),u^*(t),p_0), \\
& & p^*(t_f^*) \perp T_{\cC}(y_x^*(t_f^*)), \\
& & \disp u^*(t)=\arg\min\limits_{v\in U} H(y_x^*(t),p^*(t),v,p_0) \quad \mbox{for a.e. } t \in [0,t_f^*],
\end{eqnarray}
\end{subequations}
where $T_{\cC}(\xi)$ denotes the contingent cone of $\cC$ at $\xi$. 
Moreover, if the final time $t_f^*$ is not fixed and is an optimal time, then we have the additional condition:
\begin{eqnarray}
 & & H(y_x^*(t),p^*(t),u^*(t),p_0) = 0, \quad \mbox{for } t\in (0,t_f^*).
\end{eqnarray}
Two common cases are $\mathcal{C}=\{y_f\}$ with $p^*(t_f^*)$ free, and 
$\mathcal{C}=\R^d$ with $p^*(t_f^*)=0.$ \\ 

Now we assume that minimizing the Hamiltonian provides the control as a function $\gamma$ of the state and costate.
For a given value of $p(0)$, we can integrate $(y,p)$ by using the control $u = \gamma(y,p)$ on $[0,t_f]$.
We define the shooting function $S$ that maps the unknown $p(0)$ to the value of the final and transversality conditions at $(y(t_f),p(t_f))$.
Finding a zero of $S$ gives a trajectory $(y,u)$ that satisfies the necessary conditions for the problem ($P$).
This is typically done in practice by applying a (quasi-)Newton method.
\begin{Remark}
The multiplier $p_0$ could be equal to $0$. In that case, the PMP is said anormal, its solution 
$(y^*,p^*,u^*)$ corresponds to a ''singular'' extremal which does not depend on the cost function
$\ell$. Several works have been devoted to the existence (or nonexistence) of such extremal curves
\cite{Fr09,ChJeTr}. For numerics, in general we assume that $p_0\neq0$ which leads to solve the
PMP system with $p_0=1$. In the sequel, we will always assume that we are in the normal case ($p_0=1$).
\end{Remark}

\paragraph{Singular arcs.}
A singular arc occurs when minimizing the Hamiltonian fails to determine the optimal control $u^*$ on a whole 
time interval. The typical context is when $H$ is linear with respect to $u$, with an admissible set of controls 
of the form $U=[u_{low},u_{up}]$.
In this particular case, the function $(y,p,u) \longmapsto H_u(y,p,u)$ does not depend on the control variable. 
We define the switching function $\psi(y,p) = H_u(y,p,u)$ and have the following bang-bang control law:
$$
\left \lbrace
\begin{array}{ll}
\text{if}\ \psi(y,p) > 0 & \text{then}\ u^* = u_{low}\\
\text{if}\ \psi(y,p) < 0 & \text{then}\ u^* = u_{up}\\
\text{if}\ \psi(y,p) = 0 & \text{then switching or singular control}.
\end{array}
\right .
$$
A singular arc then corresponds to a time interval where the switching function $\psi$ is zero.
The usual way to obtain the singular control is to differentiate $\psi$ with respect to $t$ until the control 
explicitly appears, which leads to solving an equation of the form $\psi^{(2k)}(y,p) = 0$, see \cite{BryHo75}.
This step can be quite difficult in practice, depending on the problem.
Moreover, it is also required to make assumptions about the control structure, more precisely to fix the number 
of singular arcs. Each expected singular arc adds two shooting unknowns $(t_{entry},t_{exit})$, with the 
corresponding junction conditions $\psi(t_{entry}) = \dot \psi(t_{entry}) = 0$ or alternatively 
$\psi(t_{entry}) = \psi(t_{exit}) = 0$.
The problem studied in section \ref{problem3} presents such a singular arc.

\paragraph{State constraints.}
We consider a state variable inequality constraint \mbox{$g(y(t)) \le 0$}.
We denote by $q$ the smallest order such that $g^{(q)}$ depends explicitly on the control $u$; $q$ is called the order of the constraint $g$.
The Hamiltonian is defined with an additional term for the constraint
$$
H(y,p,u) = \ell(y,u) + <p,f(y,u)> + \mu g^{(q)}(y,u)
$$
with the sign condition
$$
\left \lbrace
\begin{array}{ll}
\mu = 0   & \textrm{if}\ g < 0\\
\mu \ge 0 & \textrm{if}\ g = 0.
\end{array}
\right .
$$
When the constraint is inactive we are in the same situation as for an unconstrained problem.
Over a constrained arc where $g(y) = 0$, we obtain the control from the equation $g^{(q)}(y,u)=0$, and $\mu$ from the equation $H_u=0$.
As in the singular arc case, we need to make assumptions concerning the control structure, namely the number of constrained arcs.
Each expected constrained arc adds two shooting unknowns $(t_{entry},t_{exit})$ with the Hamiltonian continuity as corresponding conditions.
We also have the so-called tangency condition at the entry point
$$
N(y(t_{entry})) := (g(y(t_{entry})), \ldots, g^{(q-1)}(y(t_{entry}))) = 0,
$$
with the costate discontinuity 
$$
p(t_{entry}^+) = p(t_{entry}^-) - \pi N_y(y(t_{entry}))   
$$
where $\pi \in \R^q$ is another multiplier yielding an additional shooting unknown.

\begin{Remark}
The tangency condition can also be enforced at the exit time, in this case the costate jump occurs at the exit time as well.
\end{Remark}

\subsection{Hamilton-Jacobi-Bellman approach}\label{subsecHJB}
 Consider the value function $\cT:\R^d\rightarrow \R$, which maps every initial condition 
$x\in \R^d$ to the minimal value of the problem ($P$). It is well known (see for example \cite{BCDbook} for a comprehensive introduction) that the value function $\cT$ satisfies
a Dynamic Programming Principle and
the Kru\v zkov transform of $\cT$, defined by 
$$
v(x):=1-e^{-\cT(x)},
$$
is the unique solution (in {\em viscosity} sense \cite{BCDbook}) of the following HJB equation
\begin{equation}\label{HJB}
\left\{
\begin{array}{ll}
v(x)+\sup\limits_{u\in U}\{-f(x,u)\cdot D v(x) -\ell(x,u)+(\ell(x,u)-1)v(x))\}=0 & x\in\R^d\backslash\cC \\
v(x)=0 & x\in\cC .
\end{array}
\right . 
\end{equation}
\\
Obtaining a numerical approximation of the function $v$ is a difficult task, mainly because $v$ is not always 
differentiable. Several numerical schemes have been studied in the literature. In this paper we will use a first-order semi-Lagrangian (SL) scheme \cite{F_appendix,F_games}. This choice is motivated by the fact that SL scheme seems the best one in order to approximate the gradient of the value function, this being our goal as we will see in the next section. More precisely, a first-order finite difference scheme is less accurate (but faster) because is not able to follow the characteristics along diagonal directions, his stencil being limited to the four neighbouring nodes (plus the considered node itself), see for example \cite{CF07}. Ultra-Bee scheme is accurate, but the solution is quite stair-shaped and not suitable for the approximation of the partial derivatives, see for example \cite{BMMZ06,BCZ09}.

We fix a (numerical) bounded domain $\Omega\supset\cC$ and we discretize it by means of a regular grid $G=\{x_i, i=1,\ldots,N_G\}$, where $N_G$ is the total number of nodes. We denote by $\widetilde v(x;h,k,\Omega)$ the fully discrete approximation of $v$, $h$ and $k$ being two discretization parameters (the first one can be interpreted as a time step to integrate along characteristics and the second one is the usual space step). We impose state constraint boundary conditions on $\partial\Omega$. The discrete version of (\ref{HJB}) is
\begin{equation}\label{HJBhk}
\left\{
\begin{array}{ll}
\widetilde v(x_i)=\widetilde H[\widetilde v](x_i) & x_i\in (\Omega\backslash \mathcal C)\cap G\\
\widetilde v(x_i)=0 & x_i\in\mathcal C\cap G
\end{array}
\right .
\end{equation}
where 
\begin{equation}\label{Htilde}
\widetilde H[\widetilde v](x_i):=\min\limits_{u\in U}\{\mathbb P_1\big(\widetilde v;x_i+hf(x_i,u)\big)+h\ell(x_i,u)(1-\widetilde v(x_i))\}
\end{equation}
and $\mathbb P_1\big(\widetilde v;x_i+hf(x_i,u)\big)$ denotes the value of $\widetilde v$ at the point \mbox{$x_i+hf(x_i,u)$} obtained by linear interpolation (note that the point $x_i+hf(x_i,u)$ is not in general sitting on the grid). The numerical scheme consists in iterating the fixed point sequence
\begin{equation}\label{SLscheme}
\widetilde v^{(n+1)}=\widetilde H[\widetilde v^{(n)}] \qquad n=1,2,\ldots
\end{equation}
until convergence, starting from $\widetilde v^{(0)}(x_i)=0$ on $\mathcal C$ and 1 elsewhere. To speed up the convergence we use the Fast Sweeping technique \cite{TCOZ03}. The function $\widetilde v$ is then extended to the whole space by linear interpolation. Once the function $\widetilde v$ is computed, we get easily the corresponding approximation $\widetilde{\mathcal T}$ of $\mathcal T$, and then the optimal control law  
in feedback form, see \cite{F_appendix,F_games} for details.

It is useful to note that the equation (\ref{HJB}) can also model a front (interface) propagation problem. Following this interpretation, the boundary of the target $\partial\cC$ is the front at initial time $t=0$, and the level set $\{x:\mathcal T(x)=t\}$ represents the front at any time $t>0$.

\section{Coupling HJB and PMP}

\subsection{Main connection}\label{subsec:mc}
It is known \cite{CaFra91} that for a general control problem with free end-point, if the value function
is differentiable at some point $x\in \R^d$ then it is differentiable along the optimal trajectory starting at 
$x$. Actually, \emph{the gradient of the value function is equal to the costate of Pontryagin's principle}.
In the context of minimum time problems (with target constraint), the link between the minimum time function
and Pontryagin's principle has been also investigated in several papers \cite{CFS00,CF05}, proving the same connection.

The main idea of the paper is to compute an approximation of the value function $\mathcal T$ solving the HJB equation on a rough grid, then approximate $D\mathcal T(x)$ ($x$ being the starting point) and finally use it as initial guess for $p(0)$. 
The approximation of the gradient can be obtained by standard first-order centred finite differences. 

In the case $\mathcal{T}\notin C^1(\R^d)$, it is proved in \cite{CFS00} that a connection between the two approaches 
still exists. More precisely, under some additional assumptions, we have
$$
p^*(t)\in D^+\mathcal{T}(y_x^*(t))\quad  \mbox{for } t\in [0,\cT(x)],
$$
where $D^+\mathcal{T}(x)$ is the \emph{superdifferential} of $\mathcal{T}$ at $x$ defined by
\begin{equation}\label{D+}
D^+\mathcal{T}(x):=\left\{\eta\in\R^d:\limsup_{y\rightarrow x}\frac{\mathcal{T}(y)-\mathcal{T}(x)-\eta\cdot(y-x)}{|y-x|}\leq 0\right\}.
\end{equation}
In the rest of this section we assume that $D^+\mathcal{T}(x)\neq\emptyset$. 
It is plain that we can not use finite difference approximation in order to compute $p(0)$ at the points where the value function $\mathcal{T}$ is not differentiable. Rather than that, we follow a different strategy. Let $\delta>0$ be is a small positive parameter, and $B(0,1)$ denote the unit ball in $\R^d$ centred at $0$. We first compute the vector 
\begin{equation}\label{grad=-direzmin}
\widetilde\xi^*:=\frac{\mathcal{T}(x+\delta\widetilde\zeta^*)-\mathcal{T}(x)}{\delta}\widetilde\zeta^*\,,\quad\textrm{ with }\quad
\widetilde\zeta^*:=\arg\min\limits_{\zeta\in B(0,1)}\mathcal{T}(x+\delta\zeta).
\end{equation}
Note that $\widetilde\zeta^*$ is an approximation of the direction of maximal decrease of $\mathcal T$, and $\|\widetilde\xi^*\|$ is an approximation of the directional derivative of $\mathcal T$ along the direction $\widetilde\zeta^*$. 
Since in the case $\mathcal{T}\in C^1(\R^d)$ the direction $\widetilde\xi^*$ is a first-order approximation of $D\mathcal{T}(x)$, we will use $\widetilde\xi^*$ as initial guess for the costate $p(0)$.

\begin{figure}[!h]
\begin{center}
\begin{tabular}{cc}
\scalebox{0.4}{\includegraphics{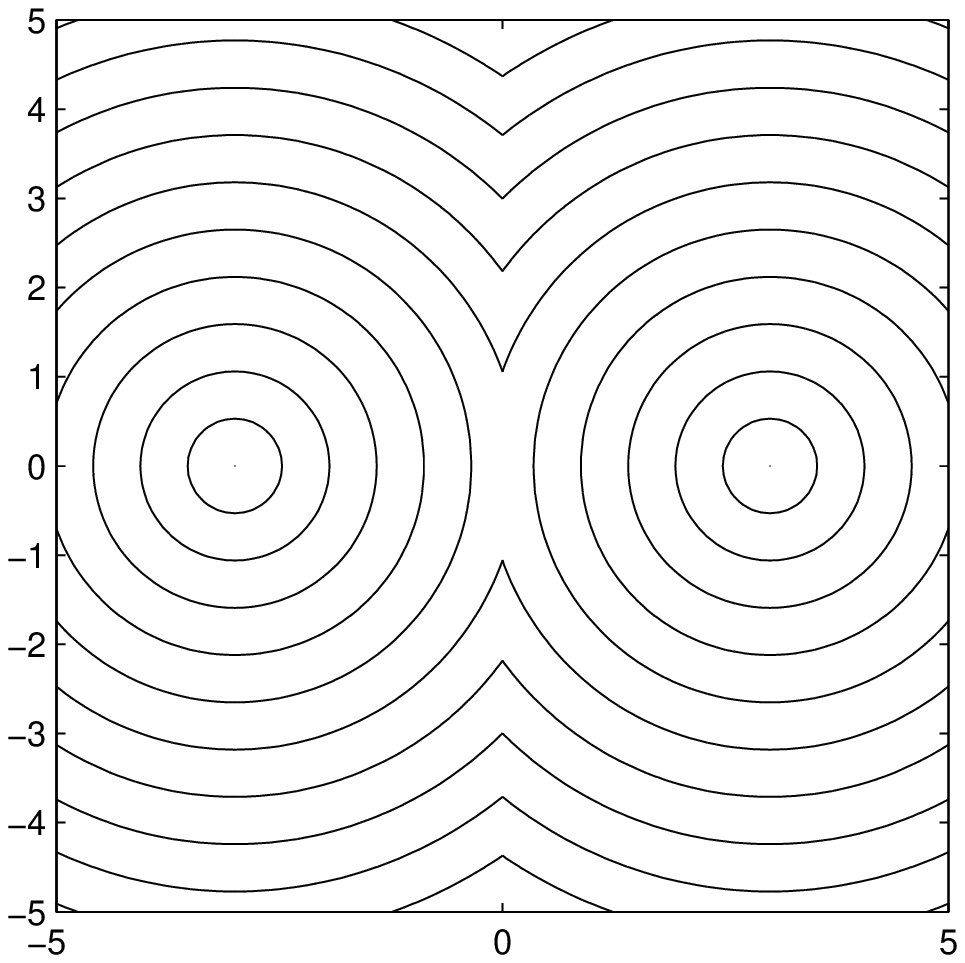}} & \hspace*{-1.5cm}
\scalebox{0.24 }{\includegraphics{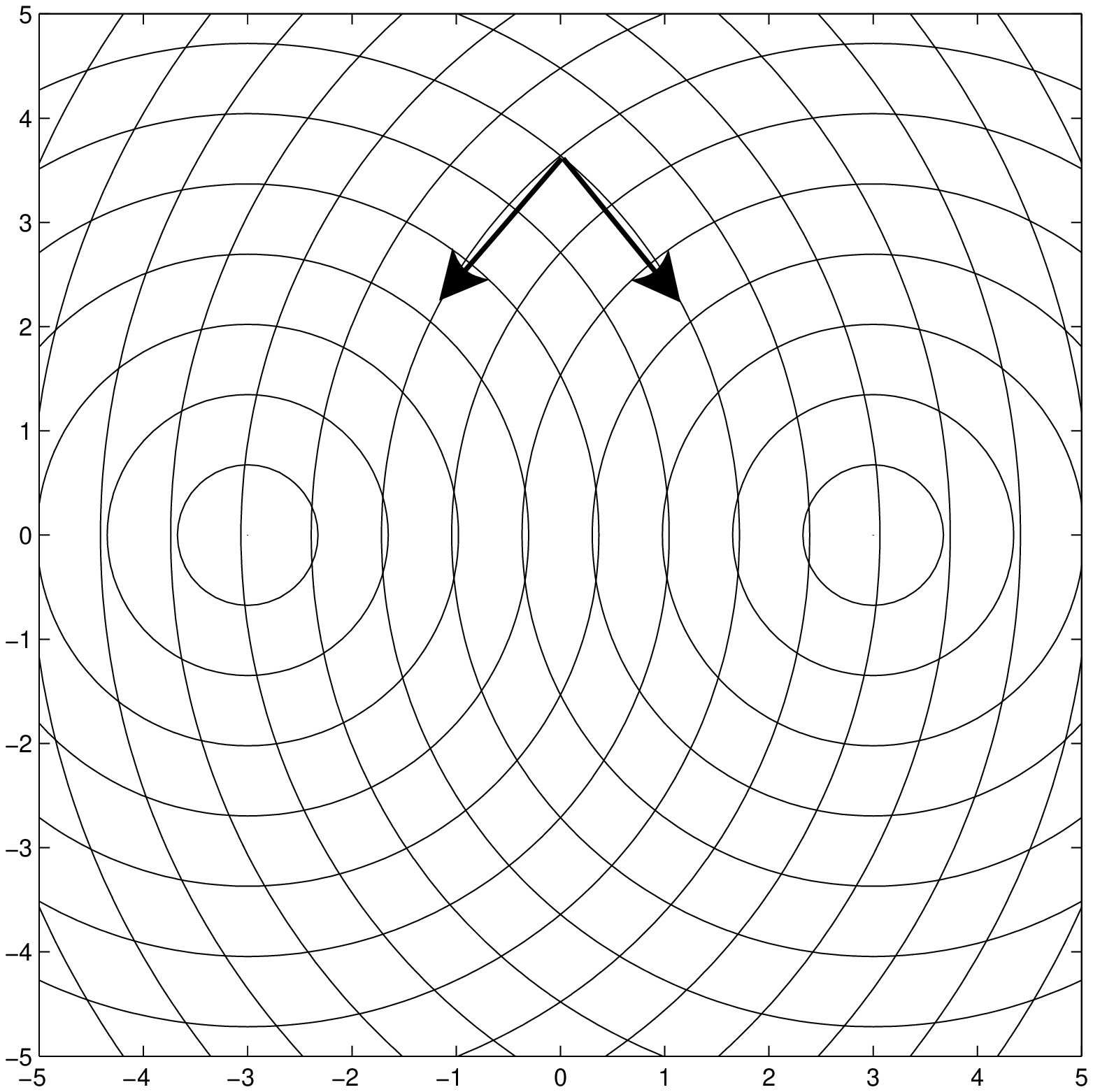}}
\end{tabular}
\caption{two crossing fronts with and without superimposition. Arrows correspond to the (two) vector(s) $\widetilde\zeta^*$}
\label{example_emi}
\end{center}
\end{figure}

Let us explain on a simple example why we choose the definition \eqref{grad=-direzmin}. Consider the case  
$$
d=2,\ \ \mathcal{C}=\{(3,0)\}\cup\{(-3,0)\},\ \ \ell \equiv 1,\ \ f=u,\ \ U=B(0,1).
$$ 
In this case the HJB equation reduces to the simple eikonal equation. 
On the line $\{x=0\}$ the function $\mathcal{T}$ is not differentiable, see the level sets in Fig. \ref{example_emi}-left. 
This line corresponds to a zone where two optimal trajectories are available, i.e. the functional $J$ has two global minima. Following the front propagation interpretation (see end of section \ref{subsecHJB}), here we have two fronts which hit each other at the line $\{x=0\}$. The viscosity solution of the HJB equation selects automatically the first arrival time so we never see the two crossing fronts, but we could 
in principle follow the propagations of the two fronts separately (see Fig. \ref{example_emi}-right), and then compute the two gradients of the two value functions. These two gradients correspond to two optimal choices for $p(0)$. By means of (\ref{grad=-direzmin}), we can approximate the two gradients without splitting the evolutions of the fronts. We first compute the two directions $\widetilde\zeta_1^*$, $\widetilde\zeta_2^*$ of maximal decrease of the function $\cT$ (see Fig. \ref{example_emi}-right), and then the ''two gradients'' $\widetilde\xi^*_1$,  $\widetilde\xi^*_2$ of $\mathcal T$.

In the present example, focusing on the point $(0,0)$, the two directions of maximal decrease
are $(-1,0)$ and $(1,0)$. It is easy to show that these two vectors coincide with the two ''extremal'' vectors in $D^+\mathcal{T}(x)$, namely the vectors $\eta$ verifying
\begin{equation}\label{D+extremal}
\limsup_{y\rightarrow x}\frac{\mathcal{T}(y)-\mathcal{T}(x)-\eta\cdot(y-x)}{|y-x|}=0.
\end{equation}
Although this relationship is not true for every function $\mathcal{T}$ such that $D^+\mathcal{T}(x)\neq \emptyset$, it is easy to see that it is true whenever the curve of non-differentiability is due to the collision of two or more fronts 
(as in Problem  1, Section \ref{problem1}). 

In section \ref{sec_ne} we will show that, beside an initial guess for $p(0)$, also other useful data can be extrapolated from the value function, and used to start the shooting method.


\subsection{Convergence of $D\mathcal T$}
Let us denote by $\widetilde D=(\widetilde D_1,\ldots,\widetilde D_d)$ the discrete gradient computed by centred finite differences with step $z>0$,
$$
\widetilde D_i\mathcal T(x):=\frac{\mathcal T(x+ze_i)-\mathcal T(x-ze_i)}{2z},\quad i=1,\ldots,d
$$
where $\{e_i\}_{i=1,\ldots,d}$ is the standard basis of $\R^d$.

Many papers (see for example \cite{BFS94,S98} in the context of differential games) investigated the convergence of the approximate value function $\widetilde v(\cdot~\!;h,k,\Omega)$ to the exact solution $v$ when the parameters $h,k$ tend to zero and $\Omega$ tends to $\R^d$. These results were quite difficult to be obtained because the function $v$ is not in general differentiable.  
To our purposes we have to go farther, proving the convergence of $\widetilde{\mathcal T}(\cdot~\!;h,k,\Omega)=-\ln(1-\widetilde v(\cdot~\!;h,k,\Omega))$ to $\cT$ and then the convergence of $\widetilde D \widetilde{\mathcal T}(\cdot~\!;h,k,\Omega)$ to $D\cT$, because the latter will be used by the PMP method as initial guess.
\\
Let us assume that $k=C_1h$ for some positive constant $C_1$. Given a generic estimate of the form
\begin{equation}\label{stima}
\|\widetilde v(\cdot~\!;h,\R^d)-v(\cdot)\|_{L^\infty(\R^d)}\leq Ch^\alpha\,, \quad C,\alpha>0
\end{equation}
we have the following
\begin{Theorem}\label{teoemi}
Assume that $\mathcal T\in C^1(\Omega)$ and there exists $\cT_{max}>0$ such that
$$
0\leq\mathcal T(x)\leq\cT_{max}~\textrm{ for all } x\in\Omega.
$$
Let us define
$$
E(x):=\|\widetilde D\widetilde{\mathcal T}(x;h,\Omega)-D\mathcal T(x)\|_\infty
$$
where $\|\ ~\|_\infty$ is the maximum norm in $\R^d$. 
\\
Then there exists $\Omega^\prime\subset\Omega$ such that
$$
\|E\|_{L^\infty(\Omega^\prime)}=O(h^\alpha /z)+O(z^2) \quad \textrm{for } h,z\rightarrow 0.
$$
\end{Theorem}
For the SL scheme we use here, an estimate of the form (\ref{stima}) in the particular case $\ell\equiv 1$ (under assumptions weaker than those used in Theorem \ref{teoemi}) can be found in \cite{S98}.
The proof of the theorem is postponed in the Appendix.

\section{Numerical experiments}\label{sec_ne}
We have tested the feasibility and relevance of combining the HJB and PMP methods on four optimal control problems.
Each of these problems highlights a particular difficulty from the control point of view.

Problem 1 (section \ref{problem1}) is a two-dimensional minimum time problem with local and 
global minima. We will see in this example that the shooting method is very sensitive with respect to the initial guess (as usual).

Problem 2 (section \ref{problem2}) is a two-dimensional controlled Van der Pol oscillator with control switchings.

Problem 3 (section \ref{problem3}) is the well-known Goddard problem with singular arcs, in the one-dimensional case (total state dimension is three).

Problem 4 (section \ref{problem4}) is another minimum time target problem in dimension four, with a first-order state constraint.

\paragraph{Details for HJB implementation.}
The algorithm is written in C++ and ran on a PC with an Intel Core 2 Duo processor at 2.00 GHz and 4GB RAM. The code is not parallelized. The reported CPU time is the time required for the whole process, which includes the computation of the value function, its gradient, the optimal trajectory and save the results on file. \\
The grid $G$ has $N_1\times \ldots \times N_d$ nodes. Grid cells have the same size $k$ in any dimension.
The set of admissible controls $U$ is discretized in $N_C$ equispaced discrete controls $\{u_j,\ j=1,\ldots,N_C\}$. 
The (fictitious) time step $h$ is variable, and chosen in such a way $h|f(x_i,u_j)|=k$ for any $x_i$ and $u_j$, so that the stencil is limited to the eight neighbouring nodes (plus the considered node itself).
The stop criterion for the fixed point iterations (\ref{SLscheme}) is 
$\|\widetilde v^{(n+1)}-\widetilde v^{(n)}\|_{L^\infty(\Omega)}<\varepsilon=10^{-5}$.

\paragraph{Details for PMP implementation.}
The algorithm is written in Fortran 90 and ran on a PC with an Intel Core 2 Duo processor at 2.33 GHz and 2GB RAM.
We used the \textsc{Shoot}\footnote{http://www.cmap.polytechnique.fr/\~{}martinon/} software which implements a shooting method with the \textsc{Hybrd} \cite{GaHiMo80} solver.
For the four problems studied we set the ODE integration method to a basic 4th-order Runge-Kutta with 100 steps.  

\subsection{Minimum time target problem}
\label{problem1}

The first example illustrates how a local solution can affect the shooting method.
We consider a simple minimum time problem in two dimensions. The goal is to reach a given position on the plane by controlling the 
angle of the velocity.
We choose the velocity in such a way  the cost functional has multiple minima.

$$
(P_1) \left \lbrace
\begin{array}{l}
\min\ t_f\\
\dot y_1(t)=c(y_1(t),y_2(t))\cos(u(t))\\
\dot y_2(t)=c(y_1(t),y_2(t))\sin(u(t))\\
U=[0,2\pi)\\  
y(0) = x = (-2.5,0)\\
y(t_f) = (3,0)
\end{array}
\right .
$$
with
$$
c(y_1,y_2)= 
\left\{
\begin{array}{ll}
1 & \textrm{if } y_2 \leq 1\\
(y_2-1)^2+1 & \textrm{if }y_2>1.
\end{array}
\right.
$$

Due to the expression of $c$, we have at least two minima.
The simplest one corresponds to a straight line trajectory ($-$) along the $y_1$ axis with $y_2=0$.
The other one has a curved trajectory ($\cap$) that takes advantage of the larger values of $y_2$. 

\subsubsection{PMP and shooting method}

We first try to solve the problem with the PMP and the shooting method.
Therefore we seek a zero of the shooting function defined by
$$
S_1: \begin{pmatrix} t_f\\ p_1(0)\\ p_2(0) \end{pmatrix} \mapsto  \begin{pmatrix}y_1(t_f)-3\\ y_2(t_f)\\ p_3(t_f)-1  \end{pmatrix}.
$$

\begin{figure}[h!]
\begin{center}
\scalebox{0.3}{\includegraphics{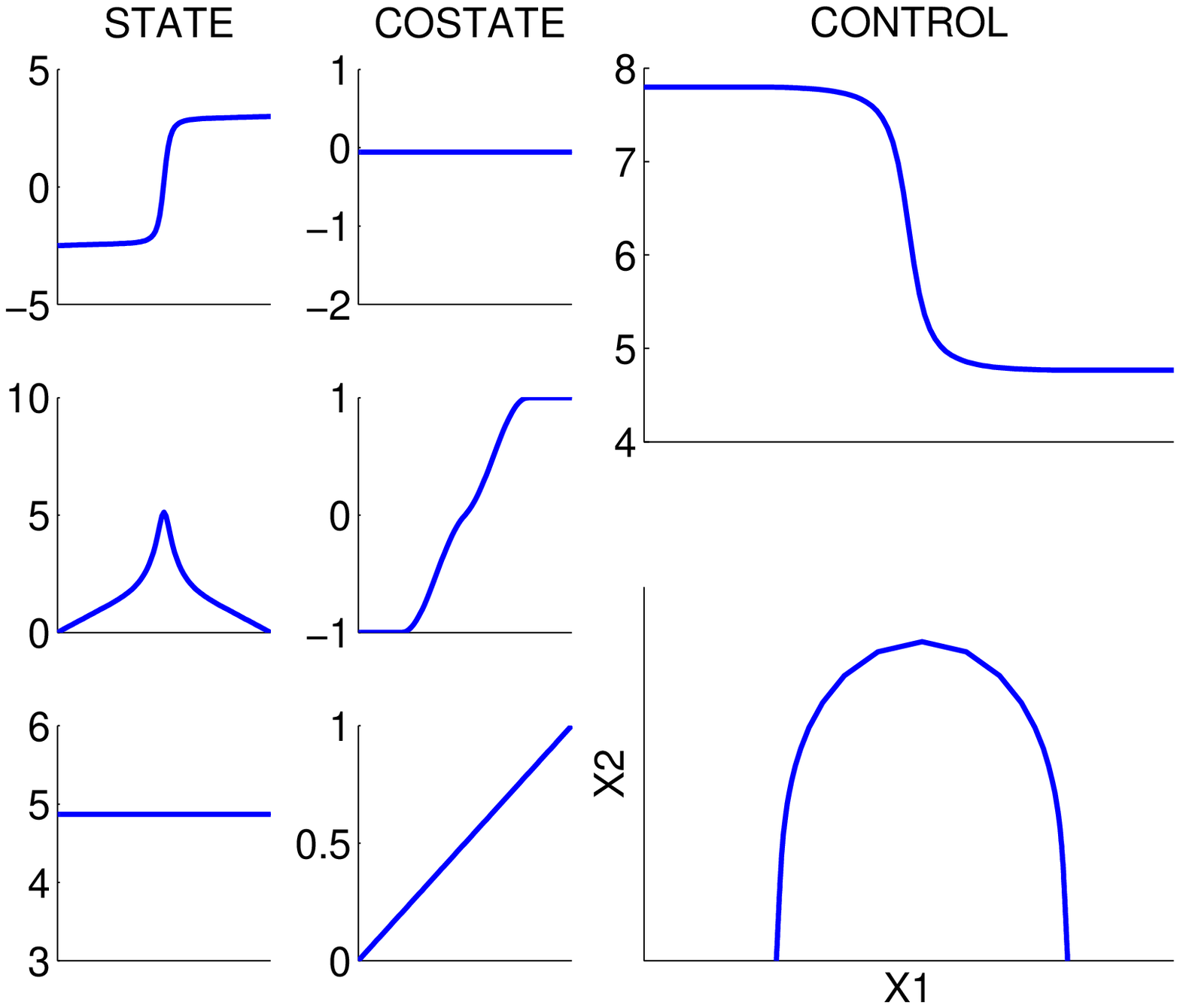}}
\hfill
\scalebox{0.3}{\includegraphics{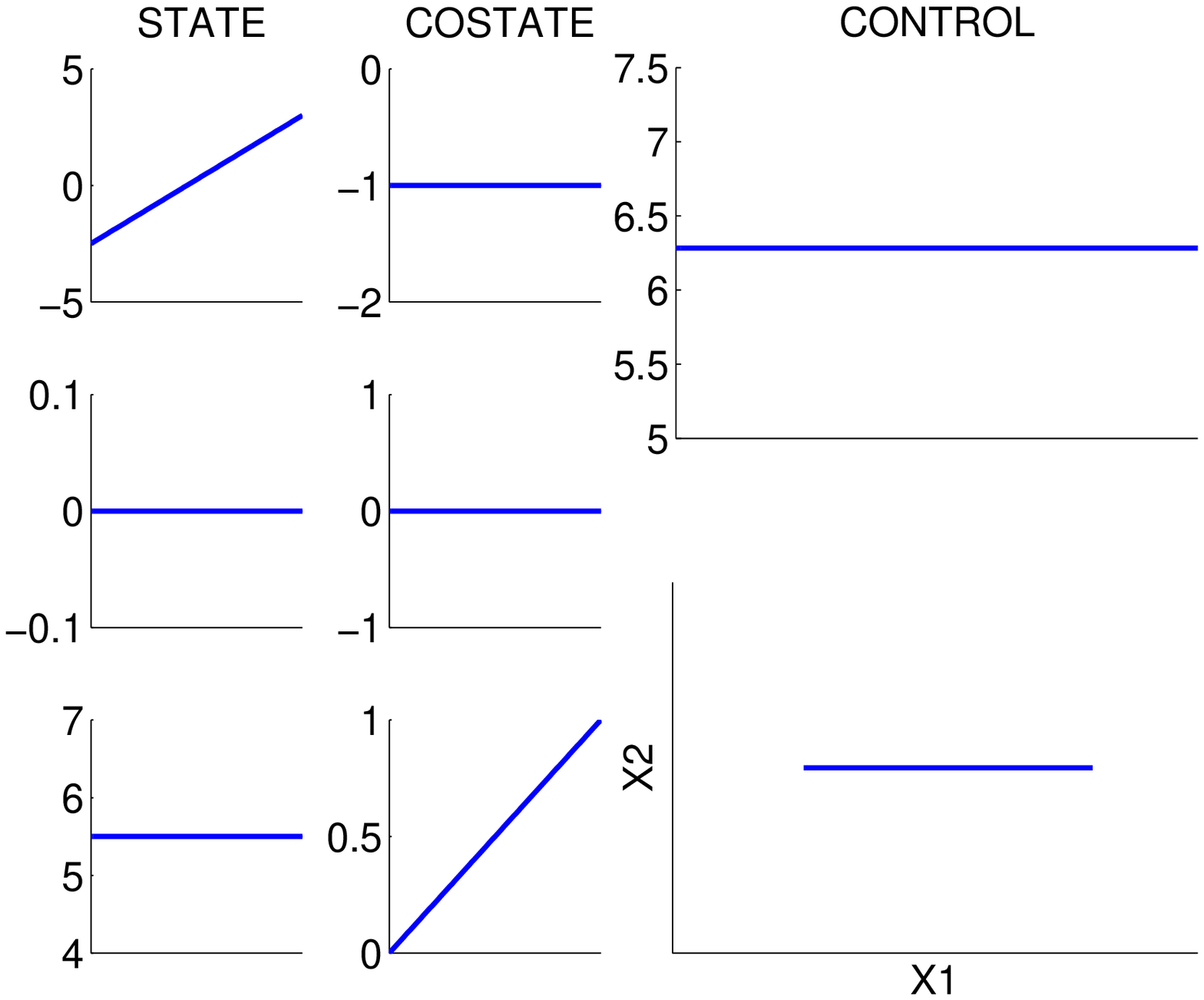}}
\caption{$(P_1)$ - Global solution (curved trajectory) and local solution (straight trajectory) found by the shooting method. Control must be intended modulo $2\pi$.} 
\label{fig:eik-glo-loc}
\end{center}
\end{figure}

\paragraph{Global and local solutions.}
Depending on the initial guess, the shooting method can converge to a local or global solution (Fig. \ref{fig:eik-glo-loc}).
The most common local solution is the straight line trajectory from $x$ to $\mathcal{C}: =\{(3,0)\}$, 
with a constant control $u=0$ and a final time $T_{local}=5.5$.
The global solution has an arch shaped trajectory that benefits from the higher speed for increasing values of $y_2$, with a final time $t_f^*=T_{global}=4.868$. 

\begin{figure}[h!]
\begin{center}
\scalebox{0.3}{\includegraphics{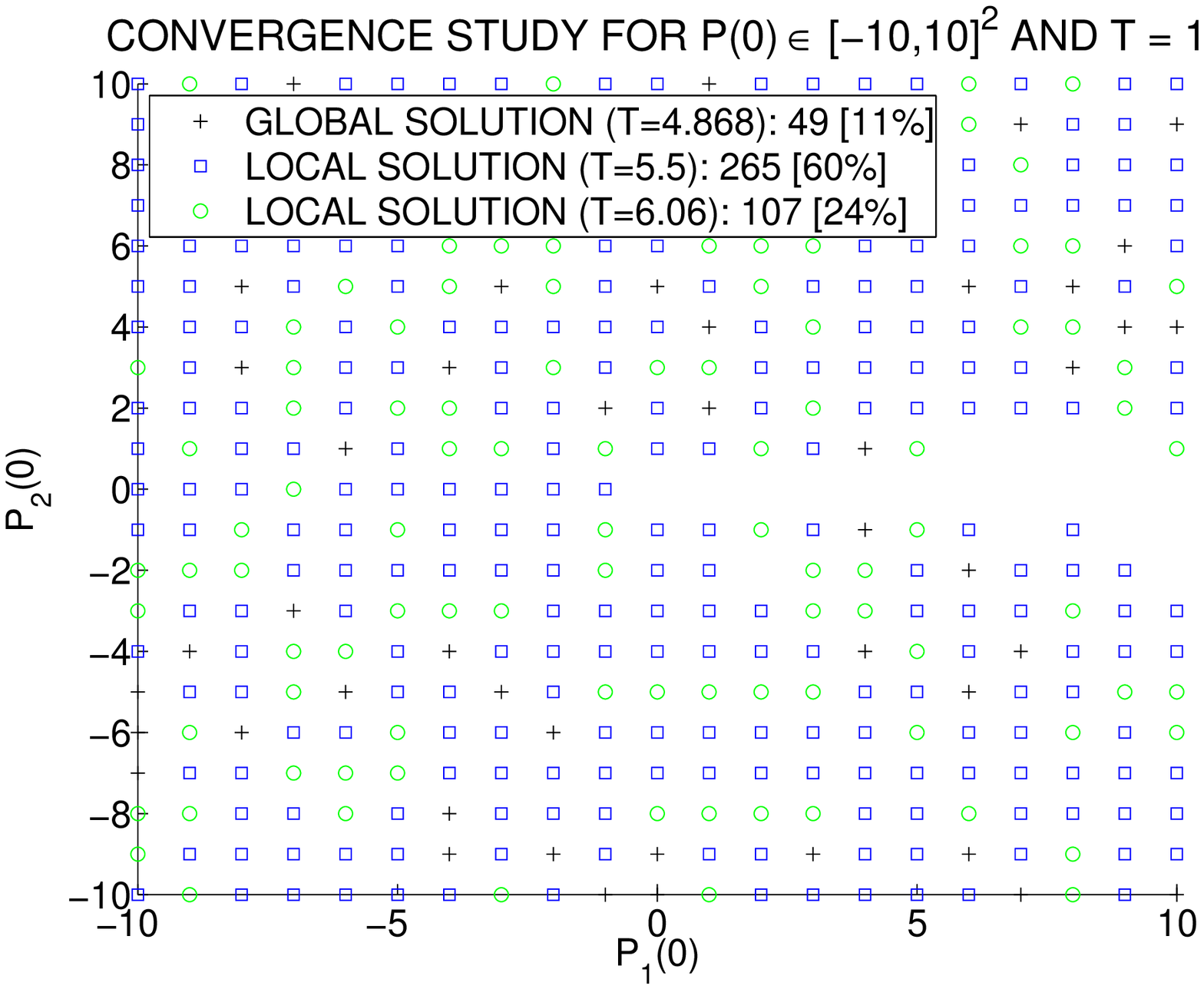}}
\hfill
\scalebox{0.3}{\includegraphics{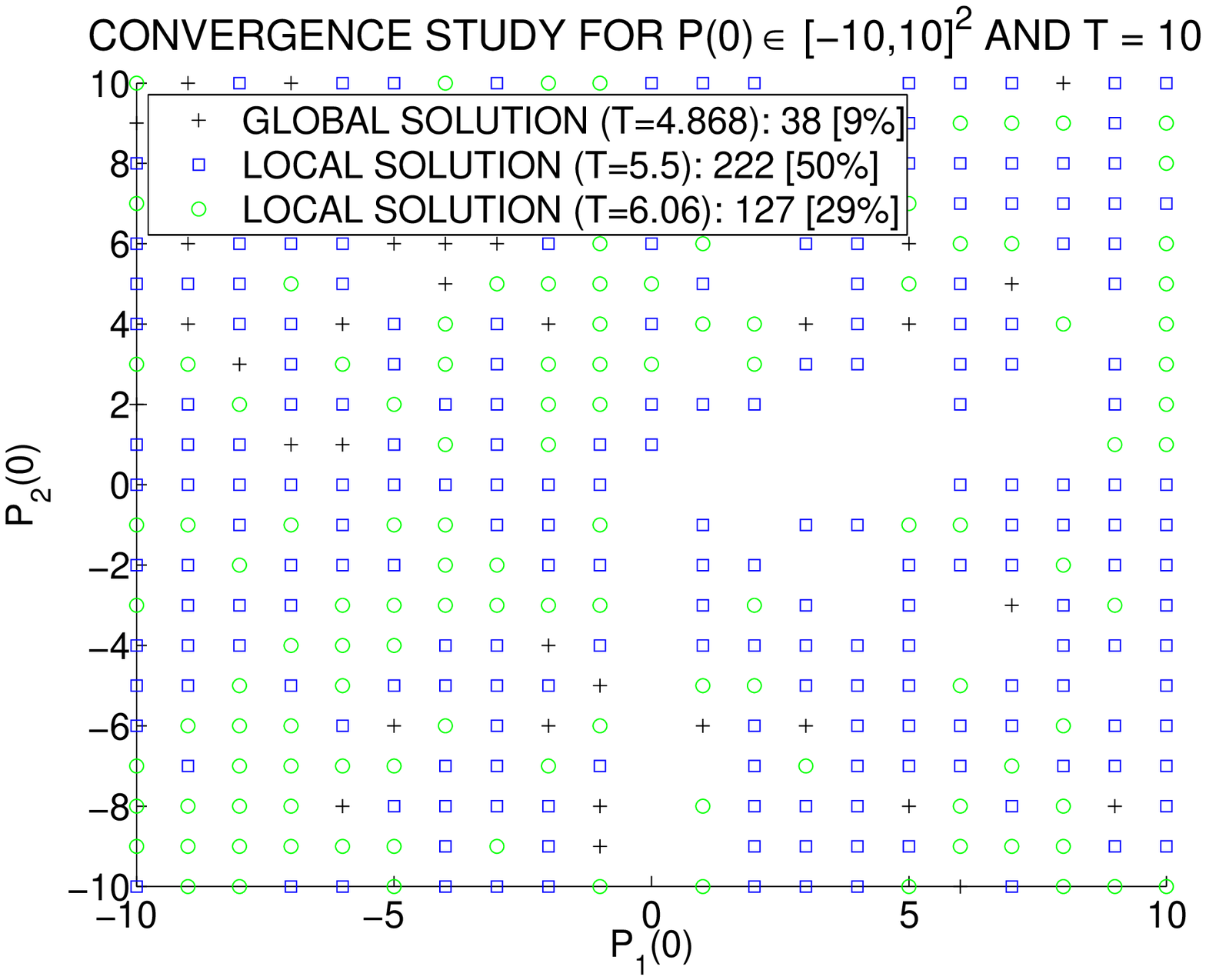}}
\caption{$(P_1)$ - Convergence from a random initialization.}
\label{fig:eik-cv}
\end{center}
\end{figure}

\paragraph{Sensitivity with respect to the starting point.}
Even for this simple problem, the shooting method is very sensitive to the starting point.
Numerical tests indicate that it converges in most cases to local solutions. 
We ran the shooting method with a batch of 441 values for $p(0)$, equally distributed in $[-10,10]^2$, and two different starting guesses for the final time (Fig. \ref{fig:eik-cv}).
We observe that for the batch with the $t_f=1$ initialization, 11\% of the shootings converge to the global solution, 60\% to the straight line local solution, and 24\% to another local solution with an even worse final time ($t_f=6.06$). Remaining 5\% does not converge at all.
For the batch with the $t_f=10$ initialization, 9\% of the shootings converge to the global solution, 50\% and 29\% to the two local solutions, and 12\% does not converge.
Obviously, just taking a random starting point is not a reliable way to find the global solution.

\subsubsection{Solving the problem with the HJB approach}

\begin{figure}[h!]
\begin{center}
\scalebox{0.55}{\includegraphics{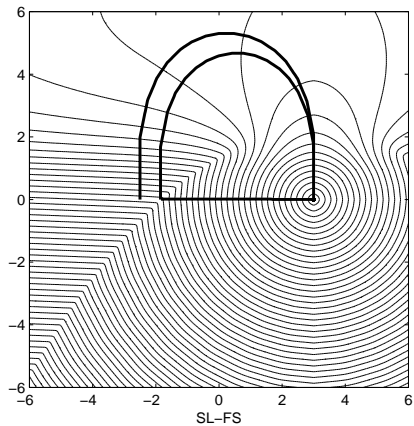}}
\caption{$(P_1)$ - Level sets of the minimum time function $\cT$, the optimal trajectory starting from
$(-2.5,0)$ and the two optimal trajectories starting from $(-1.835,0)$.}
\label{fig:eik-hjb}
\end{center}
\end{figure}

In Fig. \ref{fig:eik-hjb}, we show the level sets of the minimum time function 
$\cT$ associated to the control problem $(P_1)$. The numerical domain is $\Omega=[-6,6]^2$. 
As it can be seen, the minimum time function is not differentiable everywhere. The curve of the discontinuity of the gradient represents here the set of the initial points associated to two optimal trajectories. The superdifferential $D^+ \cT$ at the points of non-differentiability is non empty. Notice that here the minimal time function remains differentiable along each trajectory.
We will see in section \ref{problem2} a different type of non-differentiability for the value function.

\subsubsection{Coupling the HJB and PMP approaches}

We now use the data provided by the HJB approach to obtain a starting point close to the global solution.
The HJB solution provides not only an approximation of the costate $p(0)$, but also an estimate of the optimal final time $t_f^*$. 
In Table \ref{tableT1_1}, we summarize the results obtained by solving the HJB equation on several grids. As we can see, the outcome is not very sensitive with respect to the discretization parameters. This means that the choice of a very rough grid can be sufficient to obtain a good initial guess for PMP.
In fact, the shooting method immediately converges to the global solution when using the starting point 
obtained by the HJB method on the coarsest grid (Table \ref{tab:eik}).

\begin{table}[!h] 
\begin{center}
\begin{tabular}{|c|c|c|c|c|}
\hline $N_1\times N_2$ & $N_C$ & $t_f^*$ & $p(0)$  & CPU time (sec)\\
\hline\hline
$25\times 25$ & 16 & 4.895 & (-0.049,~-1.000) & 0.08  \\ \hline
$50\times 50$ & 16 & 4.895 & (-0.048,~-1.000) & 0.37 \\ \hline
$200\times 200$ & 32 & 4.878 & (-0.051,~-1.000) & 20.25 \\ \hline
\end{tabular}
\end{center}
\caption{$(P_1)$ - HJB approach: minimal time and initial costate associated to the optimal trajectory.}
\label{tableT1_1}
\end{table}

\begin{table}[h!]
\begin{center}
\begin{tabular}{|l|c|c|}
\hline
 & $t_f^*$ & $p(0)$ \\ \hline
Initialization by HJB & $4.89$     & $(-0.05, -1)$ \\
\hline
Solution by PMP       & $4.868$  & $(-5.552\times 10^{-2},-9.985\times 10^{-1})$\\
\hline
\end{tabular}
\caption{$(P_1)$ - Initialization by HJB and solution by PMP.}
\label{tab:eik}
\end{center}
\end{table}

\begin{figure}[h!]
\begin{center}
\scalebox{0.4}{\includegraphics{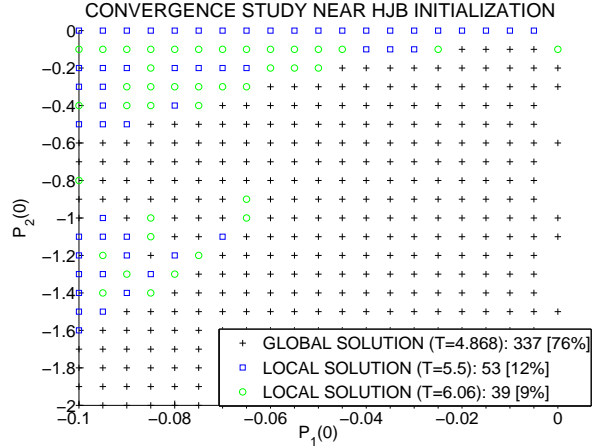}}
\caption{$(P_1)$ - Convergence to the global solution is much easier near the HJB initialization.}
\label{fig:eik-cvnearhjb}
\end{center}
\end{figure}

We can check that the convergence of the shooting method is much easier in a neighbourhood of the HJB initialization. 
We test again a batch of 441 values for $p(0)$, equally distributed in $[-0.1,0]\times[-2,0]$, which corresponds to a 100\% range around the HJB initialization $(-0.05, -1)$. We also set $t_f=4.89$.
This time the shooting method finds the global solution for $76\%$ of the points, and only $12\%$ and $9\%$ for the local solutions (Fig. \ref{fig:eik-cvnearhjb}).\\ 

We also consider the case of a starting point very close to the curve where the minimal time function is not
differentiable: $x=(-1.835,0)$. The HJB equation is solved on a $300\times 300$ grid, with $N_C=32$. The two optimal trajectories are shown in Fig. \ref{fig:eik-hjb}. Here the approximation of $p(0)$ (see section \ref{subsec:mc}) gives the two directions $p_1(0)=(-0.05,-1.00)$ and $p_2(0)=(-0.99,0.00)$.
Using these two values to initialize the shooting method, we obtain the two distinct solutions with the ``cap'' and ``straight'' trajectories (Table \ref{tab:eik-nondiff}). 
\begin{table}[h!]
\begin{center}
{
\begin{tabular}{|l|c|c|}
\hline
 & $t_f^*$ & $p(0)$ \\ \hline\hline
Initialization by HJB $(\cap)$  & $4.84$   & $(-0.05,-1)$\\ \hline
Solution by PMP $(\cap)$        & $4.8246$ & $(-7.67\times 10^{-2},-9.97\times 10^{-1})$\\
\hline\hline
Initialization by HJB $(-)$   & $4.84$   & $(-0.99,0)$\\ \hline
Solution by PMP $(-)$         & $4.835$  & $(-1,-6.2137\times 10^{-16})$\\
\hline
\end{tabular}
}
\end{center}
\caption{$(P_1)$ - Initialization by HJB and solution by PMP (two global solutions).}
\label{tab:eik-nondiff}
\end{table}


\subsection{Van der Pol oscillator}
\label{problem2}
The second test problem is a controlled Van der Pol oscillator. Here we want to reach the steady state $(y_1,y_2)=(0,0)$ in minimum time.
It is well known that the optimal trajectories, for this problem, are associated to bang-bang control variables.
$$
(P_2) \left \lbrace
\begin{array}{l}
\min\ t_f\\
\dot y_1(t) = y_2(t)\\
\dot y_2(t) = -y_1(t) + y_2(t) (1 - y_1(t)^2) + u(t)\\
U=[-1,1]\\
y(0) = x=(1,-0.8)\\
y(t_f) = (0,0)
\end{array}
\right .
$$

\subsubsection{PMP and shooting method}
Here, the Hamiltonian is linear with respect to $u$, therefore we have a bang-bang control with the switching function $\psi(y,p) = H_u(y,p,u) = p_2$. The shooting function is defined by
$$
S_2: \begin{pmatrix} t_f\\ p_1(0)\\ p_2(0) \end{pmatrix} \mapsto  \begin{pmatrix}y_1(t_f)\\ y_2(t_f)\\ p_3(t_f)-1  \end{pmatrix}.
$$
We test the shooting method with the same initial points as for problem $(P_1)$.
The convergence results are even worse in this case: for the $t_f=1$ initialization, only 9\% of the shootings converge to the global solution, and 0.5\% for the $t_f=10$ initialization.
  
\subsubsection{Solving the problem with the HJB approach}
\begin{figure}[h!]
\begin{center}
\scalebox{0.5}{\includegraphics{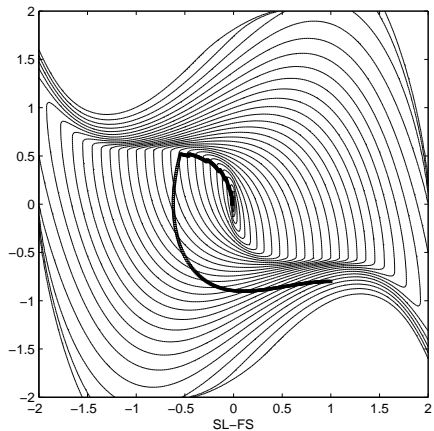}}
\caption{$(P_2)$ - Level sets of function $\cT$ and the optimal trajectory starting from $(1,-0.8)$.}
\label{fig:van-hjb}
\end{center}
\end{figure}

Here we use the HJB approach to compute the minimal time function. 
In Fig. \ref{fig:van-hjb}, we show the level sets of the solution obtained in $\Omega=[-2,2]^2$.
As in the previous problem, the value function is not differentiable everywhere,  
but here the curve of non-differentiability has a different nature. It can no more be seen as the curve of collision between two fronts and is not caused by the existence of multiple optimal solutions. It corresponds to the points where the control switches between $-1$ and $+1$. 
Taking such a starting points we have a solution with a constant control $u = \pm 1$ and no switches.
Finally we observe that at the points of non-differentiability the superdifferential is empty.

\subsubsection{Coupling the HJB and PMP approaches}
\begin{figure}[h!]
\begin{center}
\scalebox{0.22}{\includegraphics{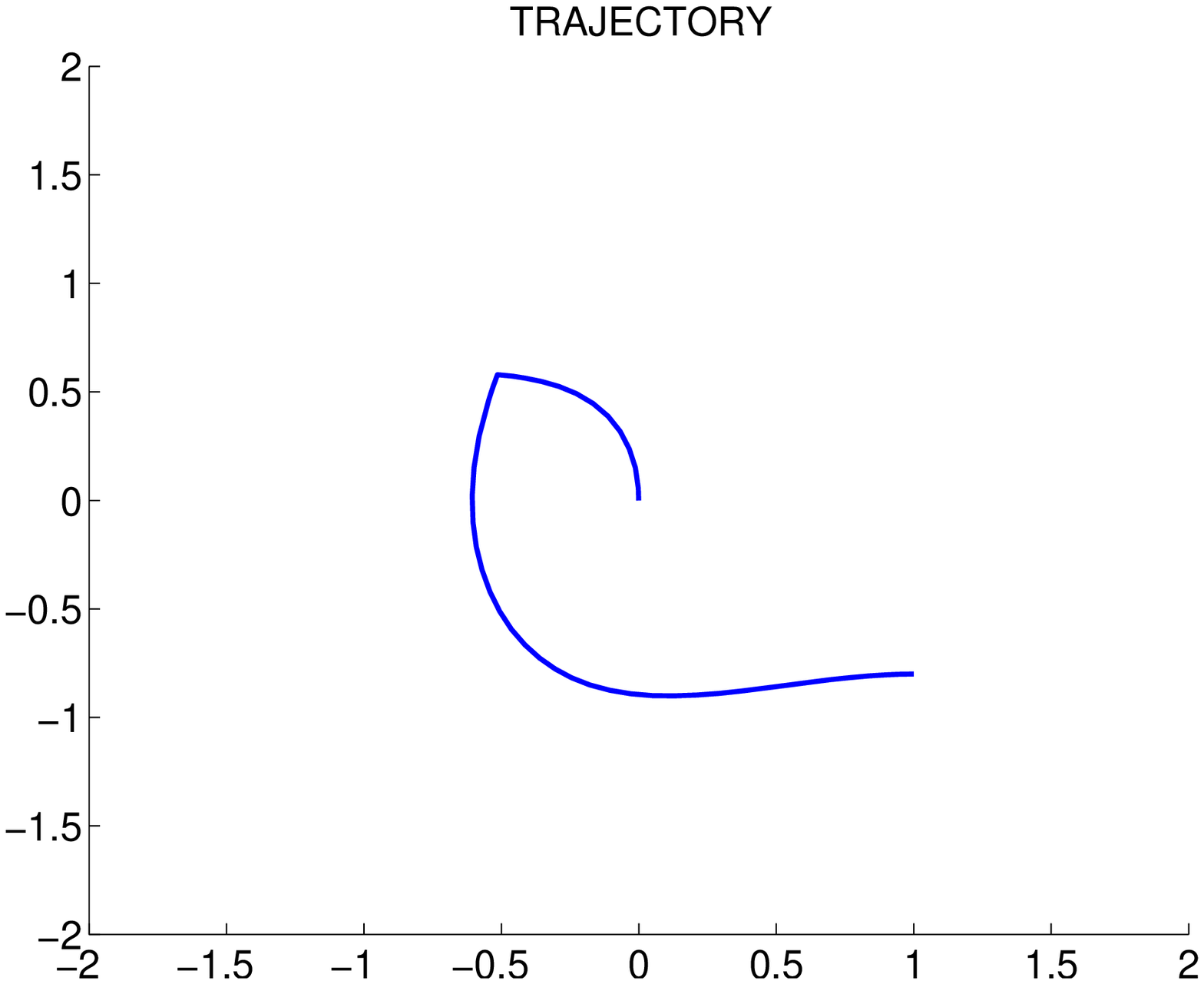}}
\scalebox{0.22}{\includegraphics{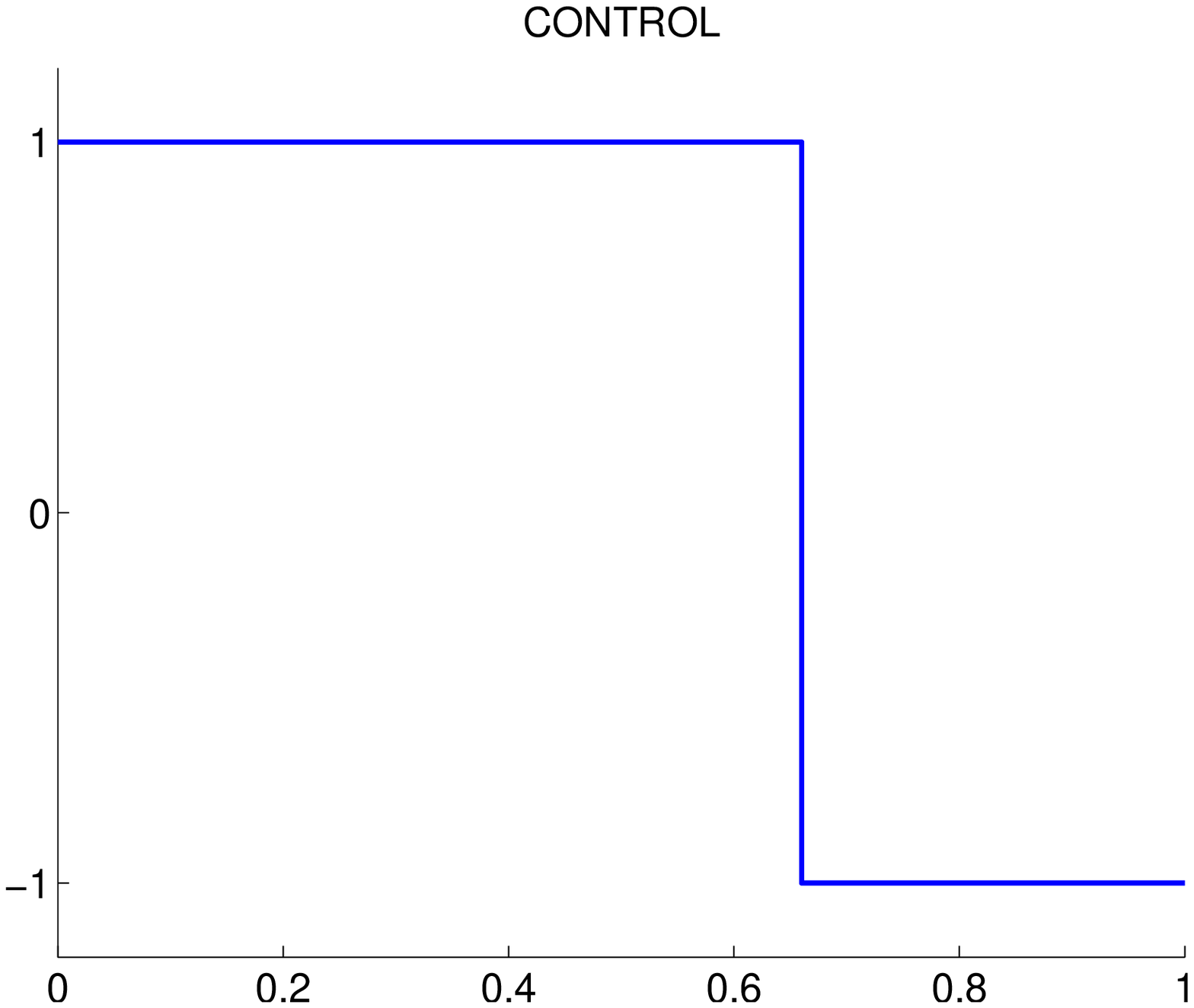}}
\scalebox{0.22}{\includegraphics{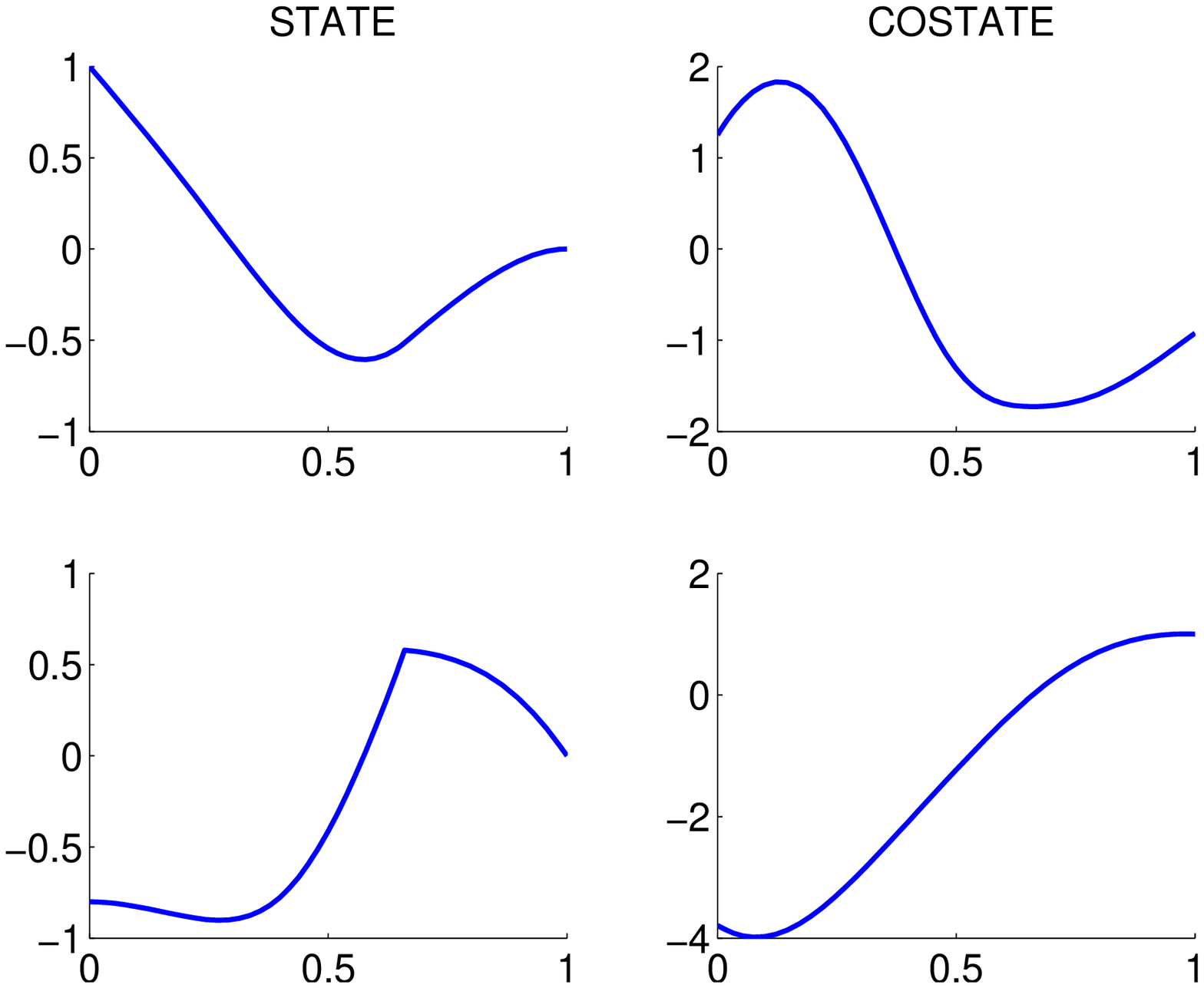}}
\caption{$(P_2)$ - Solution with one switch for the Van der Pol oscillator (shooting method).} 
\label{fig:van-pmp}
\end{center}
\end{figure}

As before, the numerical solution of the HJB equation provides an approximation of the final time $t_f$ and an initial costate $p(0)$. 
This information is used here to start the shooting algorithm.
Once again, the HJB initialization gives an immediate accurate convergence to the optimal solution, see Table \ref{tab:van} and Fig. \ref{fig:van-pmp}.
\begin{table}[h!]
\begin{center}
\begin{tabular}{|l|c|c|}
\hline
 & $t_f^*$ & $p(0)$ \\ \hline
Initialization by HJB & $4.2$     & $(1.2,-4.2)$ \\
\hline
Solution by PMP       & $3.837$  & $(1.249,-3.787)$\\
\hline
\end{tabular}
\caption{$(P_2)$ - Initialization by HJB and solution by PMP.}
\label{tab:van}
\end{center}
\end{table}
In this example, the control discontinuities hinder the convergence by testing different integration schemes for the state and costate pair $(y,p)$.
Using a fixed step integrator (4th-order Runge-Kutta) without any precaution gives a very poor convergence with a norm of $\approx 10^{-3}$ for the shooting function.
Using either a variable step integrator (\textsc{Dopri} \cite{Hairer}) or a switching detection method for the fixed step integrator  \cite{GeMa07} we get much better results ($\approx 10^{-11}$ for the shooting function norm).

\vskip0.2cm
We now test two other starting points positioned very close to the curve where the value function is not differentiable, namely $x=(1.5,-0.67)$ and $x=(1,-0.57)$.
Computation of the gradient is performed as described in section \ref{subsec:mc} in the case $\mathcal T$ is not differentiable, even if here that method is not in principle applicable due to the different nature of the non-differentiability. We observe that the shooting method finds solutions with a switch immediately after the initial time or just before the final time. 
Here the initial guesses for the costate $p(0)$ provided by the HJB method are not so close to the right ones, but they are sufficient to obtain convergence. Conversely, the minimum times given by HJB are rather close to the exact ones (Table \ref{tab:van-nondiff}).

\begin{table}[h!]
\begin{center}
\begin{tabular}{|l|c|c|}
\hline
  & $t^*_f$ & $p(0)$ \\ \hline\hline
Init. by HJB ($x=(1.5,-0.67)$)  & $3.0$    & $(1.62,-0.87)$\\\hline
Sol. by PMP  ($x=(1.5,-0.67)$)  & $2.9594$  & $(1.487,2.309\times 10^{-3})$\\
\hline\hline
Init. by HJB ($x=(1,-0.57)$)  & $2.2$     & $(1.96,-0.10)$  \\ \hline
Sol. by PMP  ($x=(1,-0.57)$)  & $2.1351$  & $(1.715,1.111\times 10^{-2})$ \\
\hline
\end{tabular}
\caption{$(P_2)$ - Initializations by HJB and solutions by PMP, non-differentiable case.}
\label{tab:van-nondiff}
\end{center}
\end{table}
%
%
%
%
%
\subsection{Goddard problem}
\label{problem3}
The third example is the well-known Goddard problem (see for instance \cite{Goddard,Mau76,Obe90,TsiKel91,SeyCli93,MaBoTre08}), to illustrate the case of singular arcs.
This problem models the ascent of a rocket through the atmosphere. We restrict ourselves to vertical (monodimensional) trajectories:
the state variables are the altitude, velocity and mass of the rocket during the flight, so we have $d=3$.
The rocket is subject to gravity, thrust and drag forces.
The final time is free, and the objective is to reach a certain altitude with a minimal fuel consumption.
$$
(P_3) \left \lbrace
\begin{array}{l}
\min ~ \int_0^{t_f} b T_{max} u\\
\dot r = v\\
\dot v = - \frac{1}{r^2} + \frac{1}{m} (T_{max} u - D(r,v))\\
\dot m = - b T_{max} u\\
U=[0,1]\\
r(0) = 1, \ v(0)=0, m(0)=1,\\
r(t_f) \ge 1.01 
\end{array}
\right .
$$
with the parameters used for instance in \cite{Obe90}: 
$b=2$, $T_{max}=3.5$ and drag $D(r,v) = 310 v^2 e^{-500(r-1)}$.\\

\subsubsection{PMP and shooting method}
As for $(P_2)$, the Hamiltonian is linear with respect to $u$, and we have a bang-bang control with possible switchings or singular arcs.
The switching function is $\psi(y,p) = H_u(y,p,u) = T_{max} ( (1-p_m)b + \frac{p_v}{m})$, and the singular control can be obtained by formally solving $\ddot{\psi}=0$.
The main difficulty, however, is to determine the structure of the optimal control, namely the number and approximate location of singular arcs.
The HJB approach is able to provide such information, in addition to the initial costate $p(0)$ and optimal time $t^*_f$.
Assuming for instance one interior singular arc, the shooting function is defined by
$$S_3: \begin{pmatrix} t_f, p_1(0), p_2(0), p_3(0) \\ t_{entry}\\ t_{exit}\end{pmatrix} \mapsto \begin{pmatrix}r(t_f)-1.01, p_2(t_f), p_3(t_f), p_4(t_f)\\ \psi(y(t_{entry}),p(t_{entry}))\\ \dot \psi(y(t_{entry}),p(t_{entry})) \end{pmatrix}.$$

\subsubsection{Solving the problem with the HJB approach}
\begin{figure}[h]
\begin{center}
\scalebox{0.45}{\includegraphics{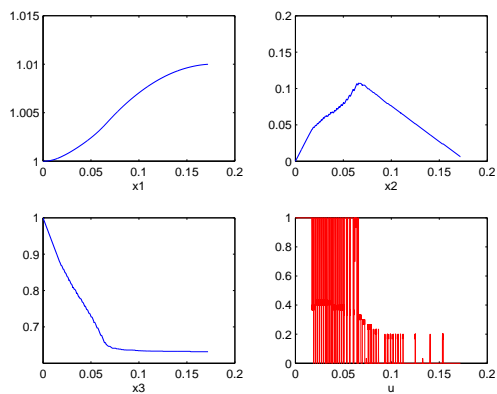}}
\caption{$(P_3)$ - Goddard problem, solution by HJB approach (first line: altitude and velocity. Second line: mass and control).}
\label{fig:god-hjb}
\end{center}
\end{figure}

The Goddard problem is also hard to solve with  the HJB approach, especially because the computation of the value function needs a huge number of iterations to converge and the solution is quite sensible to the choice of the numerical box $\Omega$ in which the value function is computed. In Fig. \ref{fig:god-hjb}, 
we show the optimal trajectory and the optimal control computed by HJB in $\Omega=[0.998,1.012]\times[-0.02,0.18]\times[0.1,1.8]$.  
As we can see, the HJB approach does not give a good approximation of the optimal control (vertical lines correspond to strong oscillations of the solution). 

\subsubsection{Coupling the HJB and PMP approaches}
\begin{figure}[h!]
\begin{center}
\scalebox{0.4}{\includegraphics{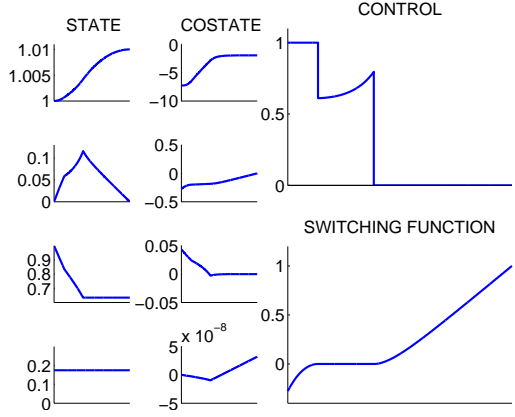}}
\caption{$(P_3)$ - Goddard problem, solution by PMP method.}
\label{fig:god-pmp} 
\end{center}
\end{figure}

As for problems $(P_1)$ and $(P_2)$, the HJB solution provides an estimate of the final time $t_f^*$ and initial costate $p(0)$.
Moreover, an examination of the HJB solution gives a good idea of the structure of the optimal control and the optimal trajectory: the change of slope on the velocity clearly visible in Fig. \ref{fig:god-hjb} indicates an interior singular arc at $(t_{entry},t_{exit}) \approx (0.02,0.06)$. The same information can be deduced by the optimal control, strongly oscillating in the same time interval.
Initializing the shooting method by means of these rough guesses, once again we obtain a quick convergence to the correct solution with the expected singular arc (Table \ref{tab:god} and Fig. \ref{fig:god-pmp}).

\begin{table}[h!]
\begin{center}
\begin{tabular}{|l|c|c|c|}
\hline
 & $t_f^*$ & $(t_{entry},t_{exit})$ & $p(0)$ \\ \hline
Init. by HJB & $0.17$ & $(0.02,0.06)$ & $(-7.79,-0.31,0.04)$ \\
\hline
Sol. by PMP       & $0.1741$  & $(0.02351,0.06685)$ & $(-7.275, -0.2773, 0.04382)$\\
\hline
\end{tabular}
\caption{$(P_3)$ - Initialization by HJB and solution by PMP.}
\label{tab:god}
\end{center}
\end{table}

\subsection{Minimum time target problem with a state constraint}\label{problem4}
This fourth example aims to illustrate the case of a state constraint, as well as a four-dimensional problem for the HJB approach.
The goal is to move a point on a plane, from an initial position to a target position, with a null initial and final velocity.
The control is the direction of acceleration, and the objective is to minimize the final time.
We add a state constraint which limits the velocity along the $x$-axis.
$$
(P_4) \left \lbrace
\begin{array}{l}
\min\ t_f\\
\dot y_1(t) = y_3(t)\\
\dot y_2(t) = y_4(t)\\
\dot y_3(t) = \cos(u(t))\\
\dot y_4(t) = \sin(u(t))\\
U=[0,2\pi)\\
y(0) = x = (-3,-4,0,0)\\
y(t_f) = (3,4,0,0)\\
y_3(t) \le 1 \quad t\in (0,t_f)
\end{array}
\right .
$$
%
Let us write the state constraint as $ g(y(t))\leq 0$, with $g$ defined by  $g(y)=y_3-1$.
The control appears explicitly in the first time derivative of $g$, so the constraint is of order 1, and we have:
$$
\dot{g}(y(t)) = \cos(u(t)), \quad g_y(y) = (0,0,1,0). 
$$
When the constraint is not active, minimizing the Hamiltonian gives the optimal control $u^*$ via 
$$(\cos(u^*),\sin(u^*)) = - \frac{(p_3,p_4)}{\sqrt{p_3^2+p_4^2}}.$$
Over a constrained arc where $g(y)=0$, the equation $\dot{g}(y,u) = 0$ and minimizing the Hamiltonian $H$ 
leads to $$u^* = - sign (p_4) \frac{\pi}{2}.$$
Then the condition $H_u = 0$ gives the value for the constraint multiplier $\mu = - p_3$.
At the entry point we have a jump condition for the costate:
$$
p(t_{entry}^+) = p(t_{entry}^-) - \pi_{entry}\ g_y, 
$$
with $\pi_{entry} \in \R$ an additional shooting unknown.
Compared to the unconstrained problem, we have three more unknowns $t_{entry},t_{exit}$ and $\pi_{entry}$.
The corresponding equations are the Hamiltonian continuity at $t_{entry}$ and $t_{exit}$ (which boils down to $p_3 = 0$), and the tangential entry condition $g(y(t_{entry})) = 0$.
The shooting function is defined by
$$
S_4: \begin{pmatrix} t_f\\ p_{1\ldots4} (0) \\ t_{entry}, t_{exit},\pi_{entry} \end{pmatrix} \mapsto \begin{pmatrix}p_5(t_f)-1\\y_{1\ldots4}(t_f)-(-3,-4,0,0)\\ p_3(t_{entry}),p_4(t_{entry}),g(y(t_{entry})) \end{pmatrix}.
$$

\begin{figure}[h!]
\begin{center}
\scalebox{0.5}{\includegraphics{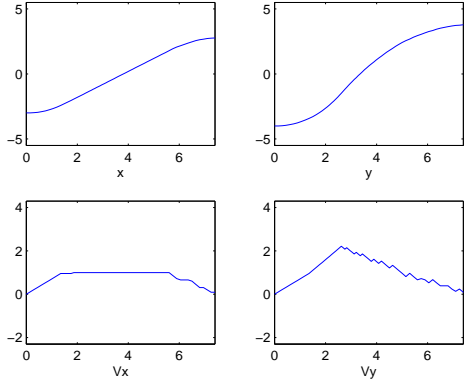}}
\caption{$(P_4)$ - Solution with a constrained arc by the HJB approach.}
\label{fig:state-hjb}
\end{center}
\end{figure}

\begin{figure}[h!]
\begin{center}
\scalebox{0.4}{\includegraphics{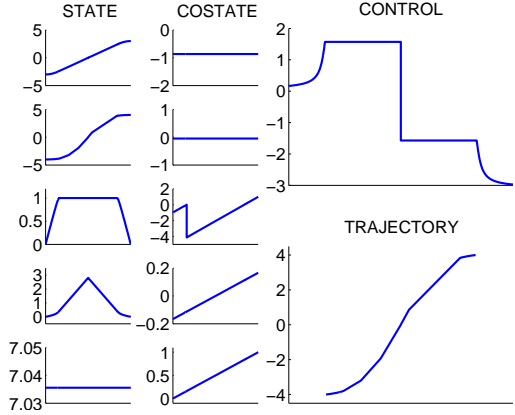}}
\caption{$(P_4)$ - Solution with a constrained arc by PMP approach.} 
\label{fig:state-pmp}
\end{center}
\end{figure}

In Fig. \ref{fig:state-hjb}, we show the numerical solution obtained by using the HJB approach in $[-5,5]^2\times[-2,4]^2$.  
In addition to the optimal final time and the initial costate, 
the HJB solution gives an estimate of the bounds for the constrained arc where $y_3=1$.
The only shooting unknown for which we were not able to obtain relevant information is the multiplier $\pi_{entry}$ for the costate jump at $t_{entry}$.
Therefore we used $\pi_{entry} = 0.1$ as a starting guess, which turned out to be sufficient for the shooting method to converge properly (Table \ref{tab:state}).
Fig. \ref{fig:state-pmp} shows the corresponding solution, much cleaner than the HJB solution but with the same structure.
We checked that the condition $\mu \ge 0$ was satisfied over the boundary arc as $p_3$ is negative, and $p_3=0$ at both entry and exit of the arc as requested by the Hamiltonian continuity conditions.
The actual value of the multiplier for the jump on $p_3$ is $\pi_{entry} = 4.1294$.

\begin{table}[h!]
\begin{center}
\begin{tabular}{|l|c|c|c|}
\hline
   & $t_f^*$ & ($t_{entry},t_{exit})$ &  $p(0)$\\
\hline
Init. by HJB & $7.5$ & $(1.35,5.6)$ & $(-0.51,-0.24,-0.89,-0.61)$   \\
\hline
Sol. by PMP  & $7.0356$ & $(1.137,5.899)$ & $(-0.867,-0.047,-0.986,-0.167)$\\
\hline
\end{tabular}
\caption{$(P_4)$ - Initialization by HJB and solution by PMP.}
\label{tab:state}
\end{center}
\end{table}

\subsection{Discretization parameters and CPU times}
In Table \ref{tab:param} we report the discretization parameters for HJB used in the four numerical tests. Note that the Van der Pol problem needs a rather fine grid to obtain a sufficient accuracy. In Table \ref{tab:cpu} we report the CPU times and the norm of the shooting function at the end of the computations. The large time for the Goddard problem is due to the huge number of iteration needed by HJB. 
\begin{table}[h!]
\begin{center}
\begin{tabular}{|l|l|l|l|l|}
\hline
 & Problem 1 & Problem 2 & Problem 3 & Problem 4\\
\hline
$N_1\times\ldots\times N_d$ & $25\times 25$ & $200\times 200$   & $20\times 20\times 20$   & $20\times 20\times 20\times 20$\\
\hline
$N_C$   & $16$    & $2$    & $20$   & $16$\\
\hline
\end{tabular}
\caption{Summary of discretization parameters for HJB}
\label{tab:param}
\end{center}
\end{table}
\begin{table}[h!]
\begin{center}
\begin{tabular}{|l|l|l|l|l|}
\hline
 & Problem 1 & Problem 2 & Problem 3 & Problem 4\\
\hline
HJB & $8\times 10^{-2}$ & $2.98$       & $211$        & $182$\\
\hline
PMP   & $3\times 10^{-3}$    & $7\times 10^{-3}$    & $3\times 10^{-2}$   & $2\times 10^{-2}$\\
\hline
$|S|$  & $2.82\times 10^{-16}$ & $8.14\times 10^{-11}$ & $1.12\times 10^{-7}$ & $6.68\times 10^{-11}$\\
\hline
\end{tabular}
\caption{Summary of CPU times (seconds) and shooting function norm}
\label{tab:cpu}
\end{center}
\end{table}
\section{Conclusions and future work}
The known relation between the gradient of the value function in the HJB approach and the costate in the PMP approach allows to use the HJB results to initialize a shooting method.
With this combined method, one can hope to benefit from the optimality of HJB and the high precision of PMP.

We have tested the combined approach on four control problems presenting some specific difficulties. 
The numerical tests also included two cases where the value function is not differentiable. 
For these four problems, the HJB approach provides useful data such as an estimate of the initial costate $p(0)$, the optimal final time $t_f^*$, and the structure of the optimal solution with respect to singular or constrained subarcs.
In each case this information allowed us to successfully initialize the shooting method.
The fact that the optimal control computed by HJB was sometimes far from the exact control did not seem to be a main issue for the shooting method initialization.
The total computational time for the combined HJB-PMP approach did not exceed four minutes, up to dimension four. 
This probably allows to run experiments in higher dimensions, like five or six. Moreover, the code for the HJB equation is easy parallelizable \cite{CFLS94}.

Even if the main limitation of the proposed method appears to be the state dimension (imposed by HJB), as we have seen
this does not mean that only simple problems can be solved.
We plan to apply this approach more specifically to trajectory
optimization for space launchers: these problems are still hard despite
having a low dimension, typically 3/4 for coplanar flight and 5/6 for a 3D
flight.
Experiments are in progress for the Ariane 5 launcher (mission toward the
GTO orbit) and for a prototype of a reusable launcher with wings (toward the LEO orbit).

\section*{Appendix}
\emph{Proof of Theorem \ref{teoemi}}.
Given the numerical domain $\Omega$ we define the set $\Omega'$ as
$$
\Omega':=\{x\in\R^d~:~\w(x;h,\Omega)\leq\min_{x'\in\partial\Omega}\w(x';h,\Omega)\}.
$$
The set $\Omega$ is the box in which the approximate solution is actually computed and $\Omega'$ represents the subset of $\Omega$ in which the solution is not affected by the fictitious boundary conditions we need to impose at $\partial\Omega$ to make computation. From the front propagation point of view, $\partial\Omega'$ represents the front at the time it touches $\partial\Omega$ for the very first time.
\\
Let us define
$v_{max}:=(1-e^{- \cT_{max}})$ and fix $x\in\Omega'$. We have
$$
\cT(x)\leq \cT_{max}<+\infty
\quad \textrm{ and } \quad
v(x)\leq v_{max}<1.
$$
By (\ref{stima}) we have 
$$
\w(x;h)\leq v(x)+C h^\alpha\leq v_{max}+C h^\alpha.
$$
Since $v_{max}<1$ there exists $h_0>0$ such that
$$
v_{max}+C h^\alpha < 1 \quad\textrm{for all } 0<h\leq h_0
$$  
then we can define
$$
\w_{max}:=v_{max}+C h_0^\alpha < 1
$$
and we have
$$
v(x)\leq v_{max}\leq \w_{max} \quad \textrm{ and } \quad \w(x;h)\leq \w_{max} \quad\textrm{for all } x\in\Omega'
\,,\quad 0<h\leq h_0.
$$
For any fixed $x\in\Omega'$, it exists $\xi_x\in[\min\{v(x),\w(x;h)\},\max\{v(x),\w(x;h)\}]$ such that
$$
|\widetilde\cT(x)-\cT(x)|=\Big|\ln\big(1- v(x)\big)-\ln\big(1- \w(x;h)\big)\Big|=\left|\frac{1}{1-~ \!\xi_x}\right||v(x)-\w(x;h)|.
$$
Since $\xi_x\leq \w_{max}$, we have
$$
|\widetilde\cT(x)-\cT(x)|\leq \frac{C h^\alpha}{1-~ \w_{max}} \quad\textrm{ for all } x\in\Omega' \textrm{ and } 0<h\leq h_0
$$
and then it exists a positive constant $C_2$ which depends by the problem's data and on $\Omega$ such that
\begin{equation}\label{estimateEmi}
\|\widetilde\cT-\cT\|_{L^\infty(\Omega')}\leq C_2 h^\alpha \quad\textrm{ for all } 0<h\leq h_0.
\end{equation}
We are now ready to recover an estimate on the gradient of the approximate solution $\widetilde\cT$. 
By (\ref{estimateEmi}) we know that, for any $i=1,\ldots,d$
$$
\widetilde\cT(x+ze_i)=\cT(x+ze_i)+E_1 \textrm{ with } |E_1|\leq C_2 h^\alpha
$$
and 
$$
\widetilde\cT(x-ze_i)=\cT(x-ze_i)+E_2 \textrm{ with }|E_2|\leq C_2 h^\alpha.
$$
So we have
$$
\widetilde D_i \widetilde\cT(x)=\frac{\cT(x+ze_i)+E_1-(\cT(x-ze_i)+E_2))}{2z}=\widetilde D_i \cT(x)+\frac{E_1-E_2}{2z}
$$ 
so that
$$
|\widetilde D_i \widetilde\cT(x)-\widetilde D_i \cT(x)|\leq \left|\frac{E_1-E_2}{2z} \right|\leq C_2\frac{ h^\alpha}{z}
$$
and then
$$
\|\widetilde D \widetilde\cT(x)-\widetilde D \cT(x)\|_\infty \leq C_2\frac{ h^\alpha}{z}.
$$
We finally obtain, for $x\in\Omega^\prime$ and $0<h\leq h_0$,
$$
\|\widetilde D \widetilde\cT(x)-D \cT(x)\|_\infty\leq \|\widetilde D \widetilde\cT(x)-\widetilde D \cT(x)\|_\infty + \|\widetilde D \cT(x)-D \cT(x)\|_\infty =
O\left(\frac{ h^\alpha}{z}\right)+O(z^2)
$$
and the conclusion follows.
$\hfill\square$



\begin{thebibliography}{50} 

\bibitem{NUDOCCCS}
B\"uskens, C.: 
NUDOCCCS –- User's manual. Universit\"at M\"unster (1996)

\bibitem{MUSCOD}
Diehl, M., Leineweber, D.B., Sch\"afer, A.A.S.: 
MUSCOD-II Users' Manual. University of Heidelberg, Germany (2001) 

\bibitem{IPOPT}
W\"achter, A., Biegler, L.T.:  
On the implementation of a primal-dual interior point filter line search algorithm for large-scale nonlinear programming, 
Mathematical Programming 106, 25--57 (2006) 

\bibitem{CBG02}
Chudej, K., B\"uskens, Ch., Graf, T.: 
Solution of a hard flight path optimization problem by different optimization codes. 
In: Breuer, M., Durst, F., Zenger, C. (eds.) High-Performance Scientific and Engineering Computing. Lecture Notes in Computational Science and Engineering 21, pp. 289--296. Springer, New York (2002)

\bibitem{MP08}
Maurer, H., Pesch, H.J.: 
Direct optimization methods for solving a complex state-constrained optimal control problem in microeconomics. 
Applied Mathematics and Computation 204, 568--579 (2008)

\bibitem{ClVi87}
Clarke, F., Vinter, R.B.:
The relationship between the maximum principle and dynamic programming.
SIAM J. Control Optim. 25, 1291--1311 (1987)

\bibitem{CaFra91}
Cannarsa, P., Frankowska, H.:
Some characterizations of optimal trajectories in control theory.
SIAM J. Control Optim. 28, 1322--1347 (1991) 

\bibitem{CFS00}
Cannarsa, P., Frankowska, H., Sinestrari, C.:
Optimality conditions and synthesis for the minimum time problem.
Set--Valued Analysis 8, 127--148 (2000)

\bibitem{CF05}
Cernea, A., Frankowska, H.:
A connection between the maximum principle and dynamic programming for constrained control problems.
SIAM J. Control Optim. 44, 673--703 (2006)

\bibitem{Pesch}
Pesch, H.J.:
A practical guide to the solution of real-life optimal control problems.
Control and Cybernetics 23, 7--60 (1994)

\bibitem{Pesch_ref}
Pesch, H.J., Plail, M.:
The maximum principle of optimal control: a history of ingenious idea and missed opportunities.
to appear in Control and Cybernetics.

\bibitem{Deu04}
Deuflhard, P.:
Newton Methods for Nonlinear Problems.
Springer Series in Computational Mathematics 35, (2004)

\bibitem{Bellman} 
Bellman, R. E.:
The theory of Dynamic Programming.
Bull. Amer. Math. Soc. 60, 503--515 (1954)

\bibitem{F_appendix}
Falcone, M.:
Numerical solution of dynamic programming equations.
Appendix A in \cite{BCDbook}.

\bibitem{F_games}
Falcone, M.:
Numerical methods for differential games based on partial differential equations.
International Game Theory Review 8, 231--272 (2006)

\bibitem{Pontryagin}
Pontryagin, L.S., Boltyanski,  V.G., Gamkrelidze, R.V. Mishtchenko,  E.F.:
The mathematical theory of optimal processes.
Wiley Interscience, New York (1962)

\bibitem{BryHo75}
Bryson, A.E., Ho, Y.-C.:
Applied optimal control.
Hemisphere Publishing, New-York (1975)

\bibitem{Pes89}
Pesch, H.J.:
Real-time computation of feedback controls for constrained optimal control problems II: a correction method based on multiple shooting.
Optimal Control, Applications and Methods 10, 147--171 (1989)

\bibitem{Fr09}
Bettiol, P., Frankowska, H.:
Normality of the maximum principle for nonconvex constrained Bolza problems.
J. Differential Equations 243, 256--269 (2007)

\bibitem{ChJeTr}
Chitour, Y., Jean, F., Tr\'{e}lat, E.:
Genericity results for singular curves.
Journal of Differential Geometry 73, 45--73 (2006)

\bibitem{BCDbook}
Bardi, M., Capuzzo Dolcetta, I.:
Optimal control and viscosity solutions of Hamilton-Jacobi-Bellman equations.
Birkh\"auser, Boston (1997)

\bibitem{CF07}
Cristiani, E., Falcone, M.: 
Fast semi-Lagrangian schemes for the Eikonal equation and applications. 
SIAM J. Numer. Anal. 45, 1979--2011 (2007)

\bibitem{BMMZ06}
Bokanowski, O., Martin, S., Munos, R., Zidani, H.:  
An anti-diffusive scheme for viability problems.
Applied Numerical Mathematics 56, 1147--1162 (2006) 

\bibitem{BCZ09}
Bokanowski, O., Cristiani, E., Zidani, H.:
An efficient data structure and accurate scheme to solve front propagation problems,
to appear in J. Sci. Comput.

\bibitem{TCOZ03}
Tsai, Y.R., Cheng, L.T., Osher, S., Zhao, H.:
Fast sweeping algorithms for a class of Hamilton-Jacobi equations.
SIAM J. Numer. Anal. 41, 673--694 (2003)

\bibitem{BFS94}
Bardi, M., Falcone, M., Soravia, P.:
Fully discrete schemes for the value function of pursuit-evasion games.
In T. Basar and A. Haurie (eds), ``Advances in Dynamic Games and Applications'', Annals of the International Society of Dynamic Games 1, 89-105 (1994)

\bibitem{S98}
Soravia, P.:
Estimates of convergence of fully discrete schemes for the Isaacs equation of pursuit-evasion differential games via maximum principle.
SIAM J. Control Optim. 36, 1--11 (1998)

\bibitem{GaHiMo80}
More, J.J., Garbow, B.S., Hillstrom, K.E.: 
User Guide for MINIPACK-1. 
Argonne National Laboratory Report ANL-80-74 (1980)

\bibitem{Hairer}
Hairer, E., N{\o}rsett, S.P., Wanner, G.:
Solving ordinary differential equations I.
Springer Series in Computational Mathematics 8, Springer-Verlag, Berlin (1993)

\bibitem{GeMa07}
Gergaud, J., Martinon. P.:
Using switching detection and variational equations for the shooting method.
Optim. Control Appl. Meth. 28, 95--116 (2007)

\bibitem{Goddard}
Goddard, R.H.:
A Method of reaching extreme altitudes.
Smithsonian Inst. Misc. Coll. 71, (1919)

\bibitem{Mau76}
Maurer, H.:
Numerical solution of singular control problems using multiple shooting techniques.
J. Optim. Theory Appl. 18, 235--257 (1976)

\bibitem{Obe90}
Oberle, H.J.:
Numerical computation of singular control functions in trajectory optimization.
Journal of Guidance, Control and Dynamics 13, 153--159 (1990)

\bibitem{TsiKel91}
Tsiotras, P., Kelley, H.J.:
Drag-law effects in the Goddard problem.
Automatica 27, 481-490 (1991)

\bibitem{SeyCli93}
Seywald, H., Cliff, E.M.:
Goddard problem in presence of a dynamic pressure limit.
Journal of Guidance, Control, and Dynamics 16, 776--781 (1993)

\bibitem{MaBoTre08}
Bonnans, F., Martinon, P., Tr\'elat, E.:
Singular arcs in the generalized Goddard's problem.
J. Optim. Theory Appl. 139, 439--461 (2008)

\bibitem{CFLS94}
Camilli, F., Falcone, M., Lanucara, P., Seghini: A.:
A domain decomposition method for Bellman equations. In: Keyes, D.E., Xu, J. (eds.), Domain Decomposition Methods in Scientific and Engineering Computing. Contemporary Mathematics 180, pp. 477--483. AMS (1994)

\end{thebibliography}
\end{document}